\documentclass{article}

\usepackage{arxiv}

\usepackage[utf8]{inputenc} 
\usepackage[T1]{fontenc}    
\usepackage{hyperref}       
\usepackage{url}            
\usepackage{booktabs}       
\usepackage{amsfonts}       
\usepackage{nicefrac}       
\usepackage{microtype}      
\usepackage{lipsum}
\usepackage{graphicx}
\usepackage{comment}
\usepackage[title]{appendix}
\usepackage{amsmath}
\usepackage[version=4]{mhchem}
\usepackage{siunitx}
\usepackage{longtable,tabularx}
\usepackage{algorithm}
\usepackage{algpseudocode}
\usepackage[backend=biber, style=numeric,sorting=none]{biblatex}
\newcommand{\bx}{\textbf{\emph{x}}}

\newcommand{\bC}{\textbf{\emph{C}}}

\newcommand{\bu}{\textbf{\emph{u}}}

\newcommand{\be}{\textbf{\emph{e}}}
\newcommand{\bp}{\textbf{\emph{p}}}

\newcommand{\bs}{\textbf{\emph{s}}}

\newcommand{\bH}{\textbf{\emph{H}}}

\title{On the Intersection of Two Conics}
\author{
 Michela Mancini \\
  Guggenheim
School of Aerospace Engineering\\
  Georgia Institute of Technology\\
  Atlanta, GA 30364 \\
  \texttt{mmancini32@gatech.edu} \\
   \And
John A. Christian \\
  Guggenheim
School of Aerospace Engineering\\
  Georgia Institute of Technology\\
 Atlanta, GA 30364\\
  \texttt{john.a.christian@gatech.edu} \\
}
\begin{document}
\maketitle

\section{Introduction}
Finding the intersection of two conics is a commonly occurring problem. For example, it occurs when identifying patterns of craters on the lunar surface \cite{Christian}, detecting the orientation of a face from a single image \cite{Kaminsky}, or estimating the attitude of a camera from 2D-to-3D point correspondences \cite{Ding}. Regardless of the application, the study of this classical problem presents a number of delightful geometric results.

In most of the cases, the intersection points are computed by finding the degenerate conic consisting of two lines passing through the common points \cite{Richter,Faucette}. Once a linear combination of the two conic matrices has been constructed, the solution of an eigenvalue problem provides four possible degenerate conics, of which only one coincides with the sought pair of lines. Then, the method proceeds by finding the intersection between one of the conics and the two lines. Other approaches make use of different methods, such as Gr\"{o}bner bases \cite{Bose1995} or geometric algebra \cite{Chomicki}.

Conic intersection, however, may be solved more intuitively with a convenient change of coordinates. In this work, we will consider two such coordinate changes. In the first approach, one of the conics is transformed into a parabola of the form \(y^2=x\), which reduces the intersection problem to finding the solution of a quartic. In the second approach, we instead use the concept of self-polar triangles---which, amazingly, reduces the conic intersection problem to the solution of a simple quadratic equation. 

\section{Mathematic preliminaries}
Since both the methods described in this work require a change of projective coordinates, we will recall here useful facts about conics and how a projective change of coordinate can be performed. 

\subsection{Conics and polarity relationship}
\begin{figure}[ht!]
    \centering \includegraphics[width=0.5\textwidth]{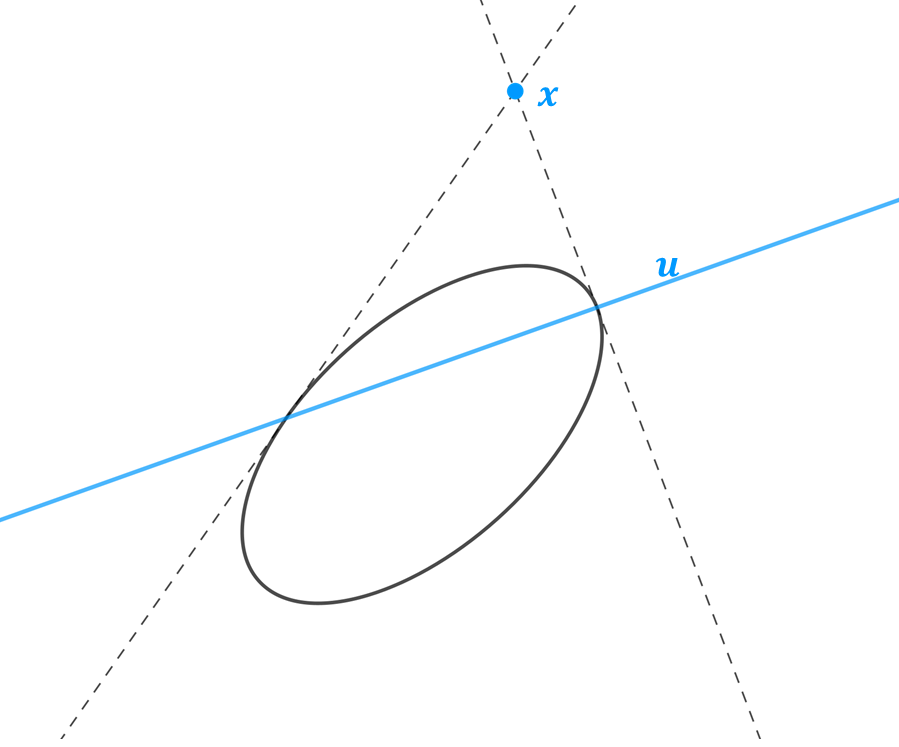}
    \caption{Relationship between a point \( {\bx}\) and its polar \( {\bu}\).}
    \label{fig:pole-polar}
\end{figure}

Consider a plane containing a conic in general position. Its equation can be written as a polynomial of degree two:
\begin{equation}\label{eq:conics}
    \,ax^2+bxy+cy^2+dx+ey+f=0
\end{equation}
Using projective coordinates, the same equation can be casted in matrix form. In fact, if we let \({\bx} \propto [x\,y\,1]^T \in \mathbb{P}^2\) be the coordinates of a point in the projective plane, the previous equation can be expressed as 
\begin{equation}
     {\bx}^T\begin{bmatrix}
         a &b/2 & d/2\\
        b/2 & c & e/2\\
        d/2 & e/2 & f\\
    \end{bmatrix} {\bx}= {\bx}^T\, \bC\,  {\bx}=0
\end{equation}
so that we can associate the conic matrix \(\bC\) of ambiguous scale
\begin{equation}
    \bC \propto \begin{bmatrix}
        a &b/2 & d/2\\
        b/2 & c & e/2\\
        d/2 & e/2 & f\\
    \end{bmatrix}
\end{equation}
The matrix \(\bC\) defines a \textit{polarity relationship} for the conic. Recalling that both points and lines can be expressed as 3-tuples in the projective plane using homogeneous coordinates (see \cite{Semple} for more details on the duality in \(\mathbb{P}^2\)), any point \( {\bx}\) on the plane can be associated to its \textit{polar line} \( {\bu}\) with respect to that conic according to
\begin{equation}
     {\bu}\propto \bC  {\bx}
\end{equation}
Similarly, if we let \(\bC^*\propto \bC^{-1}\), any line \( {\bu}\) can be associated to its polar point \( {\bx}\) using
\begin{equation}
     {\bx}\propto \bC^* {\bu}
\end{equation}
Note that \( {\bx}\) is always the pole of its polar line \( {\bu}\). The pole-polar relationship for a conic is illustrated in Fig.~\ref{fig:pole-polar}.

\subsection{Change of projective coordinates}\label{sec:coord-change}
Consider any four points \(\bp_i\) (\( i=0,\,1,\,2,\,\,3\)) in the projective plane, no three of which are collinear. Such a set of points can be used to define a new coordinate system in \(\mathbb{P}^2\). In fact, we can find a non-singular linear transformation that maps the reference points \(\be_0^T\propto[1\,0\,0]\), \(\be_1^T\propto[0\,1\,0]\), \(\be_2^T\propto[0\,0\,1]\) and the unit point \(\be_3^T\propto[1\,1\,1]\) to the four points \(\bp_0\), \(\bp_1\), \(\bp_2\) and \(\bp_3\), respectively. This mapping is represented by the homography matrix \(\bH\) such that \cite{Semple}
\begin{equation}
    \bH\be_i = \lambda_i\bp_i\quad i=0,\hdots,2
\end{equation}
where \(\lambda_i\) are unknown scale factors that need to be determined. The further condition
\begin{equation}\label{eq:lambda3}
    \bH\be_3 = \lambda_3 \bp_3
\end{equation}
fixes the scaling once the parameter \(\lambda_3\) has been specified. 
Taking advantage of our simple choices for $\be_0,\be_1,\be_2$,  the first of these two equations implies that the homography matrix \(\bH\) is 
\begin{equation}\label{eq:H}
    \bH \propto \left[\lambda_0\bp_0 \quad \lambda_1 \bp_1 \quad\lambda_2 \bp_2     \right]
\end{equation}
Selecting a scaling of \(\lambda_3 = 1\), we can re-write Eq.~\eqref{eq:lambda3} simply as \cite{Semple}
\begin{equation}
     \begin{bmatrix}
        \lambda_0 \bp_0&\lambda_1 \bp_1&\lambda_2 \bp_2 
    \end{bmatrix}\begin{bmatrix}
        1\\
        1\\
        1
    \end{bmatrix} = \bp_3
\end{equation}
which is the same as
\begin{equation}\label{eq:lambdas}
    \begin{bmatrix}
        \bp_0&\bp_1 &\bp_2
    \end{bmatrix} \begin{bmatrix}
        \lambda_0\\
        \lambda_1\\
        \lambda_2
    \end{bmatrix}=\bp_3
\end{equation}
Since \(\bp_1\), \(\bp_2\) and \(\bp_3\) are not collinear, the matrix on the left-hand side is always invertible. Thus, we can solve directly for the \(\lambda_0,\lambda_1,\lambda_2\), which allows for the computation of \(\bH\) from Eq.~\eqref{eq:H}. Given a point \( {\bx'}\) expressed in the projective coordinate system defined by \(\{\bp_i\}\), we can recover the corresponding coordinates \( {\bx}\) associated with the basis \(\{\be_i\}\) by applying
\begin{equation}
     {\bx}\propto \bH {\bx}'
\end{equation}

\section{Intersecting two conics}
In this section, we outline two different procedures that will allow us to exploit changes of projective coordinates to perform conic intersection. The first approach requires the solution of a quartic equation and is valid independently of the relative position of the two conics. The second approach requires the solution of a quadratic system, but extra care must be taken when the two conics are tangent at one single point. Pseudo-code for both the algorithms is provided in Table~\ref{alg:method2} and~\ref{alg:cap} in Appendix~\ref{sec:pseudocodes}.

\subsection{Intersecting two conics using the canonical representation}
Consider the two conics described by the matrices \(\bC_1\) and \(\bC_2\) and, without loss of generality, select any three points on the conic \(\bC_1\). Let these points be \(\bp_1\), \(\bp_2\) and \(\bp_3\). Since they need not to be real, we can select them by taking the intersection between the conic \(\bC_1\) and any two lines. For example
\begin{equation}\label{eq:points-on-conic}
    \bp_{1,2} \propto \begin{bmatrix}
        0\\
        \frac{1}{2c_1}\left(-e_1\pm\sqrt{e_1^2-4c_1f_1}\right)\\
        1
    \end{bmatrix}\quad \bp_3 \propto \begin{bmatrix}
        \frac{1}{2a_1}\left(-d_1+\sqrt{d_1^2-4a_1f_1}\right)\\
        0\\
        1
    \end{bmatrix}
\end{equation}
Different choices for the arbitrary lines would be required for $\bp_{1,2}$ if $c_1 \simeq 0$ or for $\bp_3$ if $a_1 \simeq 0$.

Now let \(\bp_0\) be the pole of the line through \(\bp_1\) and \(\bp_2\). Its coordinates can be obtained as the intersection of the two lines \(\boldsymbol{\ell}_1\propto \bC_1 \bp_1\) and \(\boldsymbol{\ell}_2\propto \bC_1\bp_2\) that are tangent to the conic at $\bp_1$ and $\bp_2$
\begin{equation}\label{eq:p0}
    \bp_0\propto \boldsymbol{\ell}_1\times \boldsymbol{\ell}_2
\end{equation}
If we let \(\bp_0\), \(\bp_1\) and \(\bp_2\) be the reference points of a new coordinate system, and \(\bp_3\) be the unit point, we can construct a matrix \(\bH\) as described in Sec.~\ref{sec:coord-change} under whose mapping the matrix \(\bC_1\) is transformed to the matrix
\begin{equation}\label{eq:c1p-parabola}
    \bC_1'\propto  \bH^T\bC_1\bH \propto \begin{bmatrix}
        2 & 0 & 0\\
        0 & 0 & -1\\
        0 & -1 & 0
    \end{bmatrix}
\end{equation}
which is the equation for a parabola.

Applying the same transformation to the matrix \(\bC_2\), we obtain the generic matrix 
\begin{equation}\label{eq:c2Generic}
    \bC_2'\propto \bH^T\bC_2\bH \propto \begin{bmatrix}
        a_2' & b_2'/2 & d_2'/2\\
        b_2'/2 & c_2' & e_2'/2\\
        d_2'/2 & e_2'/2 & f_2'
    \end{bmatrix}
\end{equation}

The intersection problem reduces now to the problem of finding the common solutions of the equations $\bx'^T \bC'_1 \bx' = 0$ and $\bx'^T \bC'_2 \bx' = 0$, which we may write as 
\begin{equation}\label{eq:parabola}
\bx'^T \bC'_1 \bx' = 0 
\quad \rightarrow \quad 
x'^2=y'
\end{equation}
\begin{equation}\label{eq:generic}
\bx'^T \bC'_2 \bx' = 0 
\quad \rightarrow \quad a_2'x'^2+b_2'x'y'+c_2'y'^2+d_2'x'+e_2'y'+f_2'=0
\end{equation}

Substituting Eq.~\eqref{eq:parabola} into Eq.~\eqref{eq:generic} and rearranging the terms we obtain a quartic equation of the form
\begin{equation}\label{eq:quartic}
    c_2'x'^4 + b_2'x'^3+ \left(a_2'+e_2'\right)x'^2 +d_2'x'+f_2'=0
\end{equation}
A quartic equation has a closed form expression (see Appendix~\ref{app:quarticsolution}). Alternatively, the solutions may be found numerically in many common software packages, such as through the \texttt{roots()} command in \texttt{MATLAB}.
Once the four solutions to Eq.~\eqref{eq:quartic} have been found, the corresponding \(y'\) are given by Eq.~\eqref{eq:parabola}, and the four pairs \((x_j',\,y_j')\) provide the coordinates of the interection points in the new coordinate system, \(\bs_j'\):
\begin{equation}\label{eq:quartic-sp}
    \bs_j'\propto \begin{bmatrix}
        x_j\\
        y_j\\
        1
    \end{bmatrix} \qquad  j=a,b,c,d
\end{equation}
From these, we can finally recover the intersection points in the initial coordinates:
\begin{equation}\label{eq:quartic-s}
    \bs_j \propto \bH \bs_j'  \qquad  j=a,b,c,d
\end{equation}
Note that if one of the points \(\bp_{0}\), \(\bp_{1}\), \(\bp_{2}\) is an intersection point, Eq.~\eqref{eq:quartic} will have degree less than four. The missing solution will be one of the reference points. 

\subsection{Intersecting two conics using a self-polar triangle.}
The solution of a quartic equation can be avoided by use of the self-polar triangle. 
Given three points in the plane, they are said to form the vertices of a self-polar triangle if the polar line of each vertex is the line passing through the remaining two vertices.

Any four distinct points on a conic may be used to form a quadrangle. By constraining these points to be on a conic, it may be shown that the quadrangle's diagonal triangle will be a self-polar triangle \cite{Semple}. While there are an infinite number of self-polar triangles for a given conic (since there are an infinite number of quadrangles from 4-tuples of distinct points on the conic), two distinct conics intersecting in four distinct points (over the complex numbers) can have only one self-polar triangle in common. In other words, there are only three points in the projective plane that form a triangle that is self-polar for both the conics \cite{Semple,Woods}. If we let these vertices be the reference points of a new projective coordinate system, the two conic equations are greatly simplified, and their intersection can be performed easily. 

\begin{figure}[ht!]
    \centering
    \includegraphics[width=0.7\textwidth]{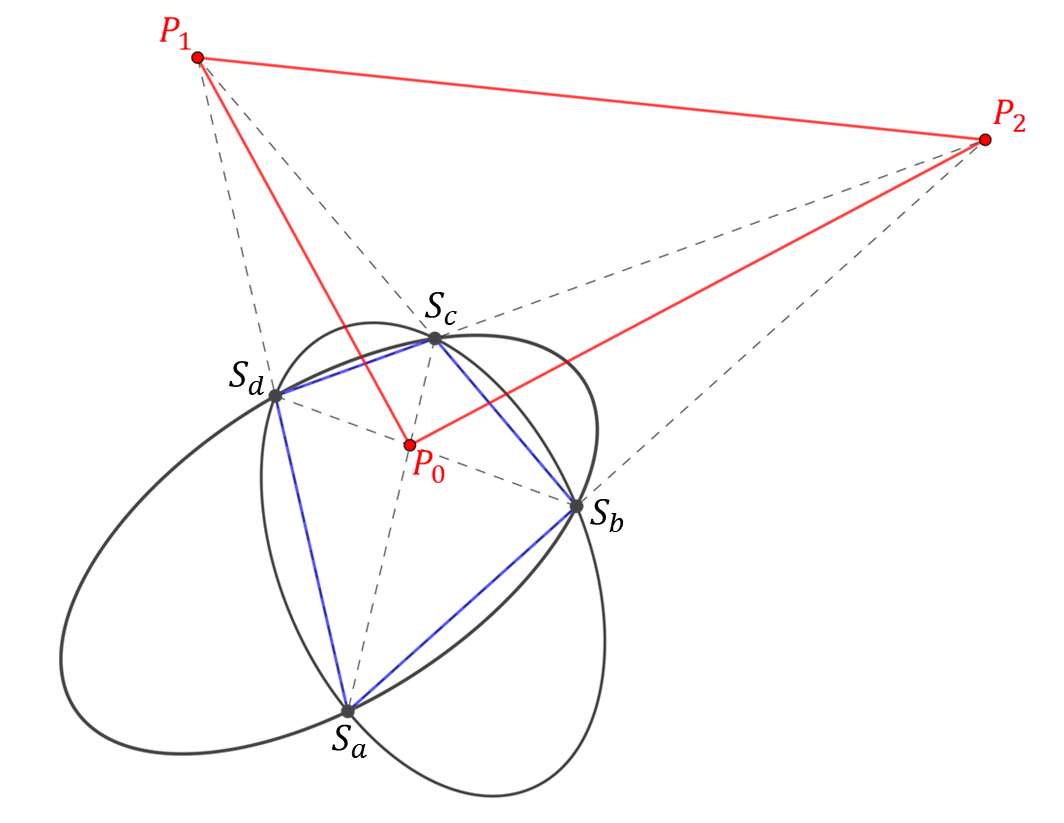}
    \caption{The points of intersection between non-consecutive edges of the quadrilateral \(S_aS_bS_cS_d\) are the vertices of the unique triangle that is self-polar for both the two conics.}
    \label{fig:self-polar}
\end{figure}

In order to visualize the common self-polar triangle, let \(S_a\), \(S_b\), \(S_c\) and \(S_d\) be the points of intersection of the two conics \(\bC_1\) and \(\bC_2\). The vertices of the unique common self-polar triangle can be found by intersecting the lines joining non-consecutive sides of the quadrilateral \(S_aS_bS_cS_d\) \cite{Semple}. This concept is illustrated in Fig.~\ref{fig:self-polar} for the case of four real intersection points. 

When the vertices of a conic's self-polar triangle are used as reference points, the conic matrix is transformed to a diagonal matrix. This can be easily seen if we consider that the reference points in the new frame have coordinates 
\begin{equation}
    \bp_0' \propto \begin{bmatrix}
        1\\
        0\\
        0
    \end{bmatrix}\quad \bp_1' \propto \begin{bmatrix}
        0\\
        1\\
        0
    \end{bmatrix}\quad\bp_2' \propto \begin{bmatrix}
        0\\
        0\\
        1
    \end{bmatrix}
\end{equation}
Since the line passing through \(\bp_l'\) and \(\bp_m'\) has coordinates \(\boldsymbol{\ell}_i \propto \bp_l'\times \bp_m'\) and \( \bp_l'\times \bp_m' \propto \bp_i'\), the polarity relationship imposes that
\begin{equation}
   \bp_i'\propto \bC' \bp_i'
\end{equation}
for all \(i\in\{0,1,2\}\), that can only happen when \(\bC'\) is diagonal. If a triangle is self-polar for two conics,  the change of coordinates will simultaneously diagonalize the conic matrices. 

In order to construct the matrix \(\bH\), we need to determine the vertices of the self-polar triangle \(\bp_0\), \(\bp_1\) and \(\bp_2\) in the original coordinate system. Given the definition of self-polar triangle, we know that the polar line of any vertex passes through the other two. Since the triangle is self-polar for the two conics, this polarity relationship must hold for both, so that we have
\begin{equation}
    \bC_1 \bp_i \propto \bC_2 \bp_i \qquad i=0,1,2
\end{equation}
In order to remove the proportionality, we can introduce the unknown scaling factor \(\mu_i\), yielding
\begin{equation}
    \bC_1 \bp_i = \mu_i\bC_2 \bp_i
\end{equation}
that is nothing more than an eigenvalue problem
\begin{equation}\label{eq:self-eig}
   \left(\bC_2^{*} \bC_1 - \mu_i \textbf{I}_{3\times 3}\right)\bp_i=\textbf{0}
\end{equation}
We can conclude that the vertices of the common self-polar triangle are the eigenvectors of the matrix \(\bC_2^{*}\bC_1\). Once that the coordinates of the reference points have been determined, we can freely choose the unit point \(\bp^T_3 = [1\,1\,1] \) so long as no three of these points are collinear. 

These values of $\bp_0, \bp_1, \bp_2,\bp_3$ may be used to construct the homography \(\bH\) using the procedure from Section~\ref{sec:coord-change}. After the transformation, the two conic matrices assume the form
\begin{equation}\label{eq:c1p}
    \bC_1'\propto \bH^T\bC_1\bH\propto \begin{bmatrix}
    a_1' & 0 & 0\\
    0 & c_1' & 0\\
    0 & 0 & f_1'
    \end{bmatrix}
\end{equation}

\begin{equation}\label{eq:c2p}
    \bC_2'\propto \bH^T\bC_2\bH\propto \begin{bmatrix}
    a_2' & 0 & 0\\
    0 & c_2' & 0\\
    0 & 0 & f_2'
    \end{bmatrix}
\end{equation}
so that the conic equations are, respectively
\begin{equation}\label{eq:self-polar1}
    a_1'x'^2+c_1'y'^2+f_1'=0
\end{equation}
and
\begin{equation}\label{eq:self-polar2}
    a_2'x'^2+c_2'y'^2+f_2'=0
\end{equation}
Their intersection is now easily found. From Eq.~\eqref{eq:self-polar1} we have
\begin{equation}\label{eq:self-polarY_sq}
    y'^2 = -f_1'/c_1'-(a_1'/c_1')x'^2
\end{equation}
that, substituted into Eq.~\eqref{eq:self-polar2}, gives
\begin{equation}
    \left(a_2'-a_1'c_2'/c_1'\right)x'^2 = f_1' c_2'/c_1'-f_2'
\end{equation}
that is 
\begin{equation}\label{eq:self-polarX_sq}
    x'^2 = {\frac{f_1'c_2'-f_2'c_1'}{a_2'c_1'-a_1'c_2'}}
\end{equation}
Using Eq.~\eqref{eq:self-polarY_sq} we can similarly write down an expression for $y'^2$ 
\begin{equation}\label{eq:self-polarY_sq_v2}
    y'^2 =\frac{f_2' a_1' - f_1' a_2'}{a_2' c_1'- a_1' c_2'}
\end{equation}
The four points of intersection (in the transformed space) are then found as
\begin{equation}\label{eq:self-s_p}
    \bs_{a,b}' \propto \begin{bmatrix}
        \pm \sqrt{x'^2}\\
        \pm \sqrt{y'^2}\\
        1
    \end{bmatrix}\qquad \bs'_{c,d}\propto \begin{bmatrix}
        \pm \sqrt{x'^2}\\
        \mp \sqrt{y'^2}\\
        1
    \end{bmatrix}
\end{equation}
and we may convert them to the original coordinate system through \begin{equation}\label{eq:self-s}
    \bs_j \propto \bH \bs_j'  \qquad  j=a,b,c,d
    \end{equation}

Two conics do not always have a unique common self-polar triangle \cite{Woods}. One exception is given by two conics tangent in two points. This case doesn't really represent a problem for our approach, as the number of common self-polar triangles is infinite and the algorithm continues to be valid. The other exception is represented by the case of two conic being tangent at one point (the other two intersections being either complex conjugates or real). In this case, they do not have any common self-polar triangle, as two of its vertices coincide, making the matrix on the left-hand side of Eq.~\eqref{eq:lambdas} singular. When this happens, however, the two coinciding vertices lie at the double intersection point, and this information may be used to determine the remaining points easily, as we will see now. 

Assume that during the construction of the self-polar triangle we encounter the situation where two vertices coincide. In this case, we can set \(\bp_1\) as one of those vertices, and re-define the other two reference points. In particular, we can choose \(\bp_2\) and \(\bp_3\) on the conic, for example as in Eq.~\eqref{eq:points-on-conic}, and construct \(\bp_0\) as in Eq.~\eqref{eq:p0}. With this choice, the first matrix will be transformed to the matrix of Eq.~\eqref{eq:c1p-parabola}, while the second matrix will transform such that \(b_2'=c_2'=0\) in Eq.~\eqref{eq:c2Generic}.The remaining two intersection points are obtained by solving the quadratic equation
\begin{equation}\label{eq:quadratic}
    (a_2'+e_2')x'^2+d_2'x'+f_2'=0
\end{equation}
that give
\begin{equation}\label{eq:quadratic-sol}
    x'=\frac{-d_2'\pm\sqrt{d_2'^2-4f_2'(a_2'+e_2')}}{2\left(a_2'+e_2'\right)}
\end{equation}

\section{Numerical example}
In this section we will show two numerical examples of application. In the first one, two conics intersect in four distinct real points. In the second case, one ellipse and a hyperbola are tangent at one point, with the remaining two intersection points being complex-conjugates.

\subsection{Four distinct intersection points}
\begin{figure}[ht!]
    \centering    \includegraphics[width=0.5\textwidth]{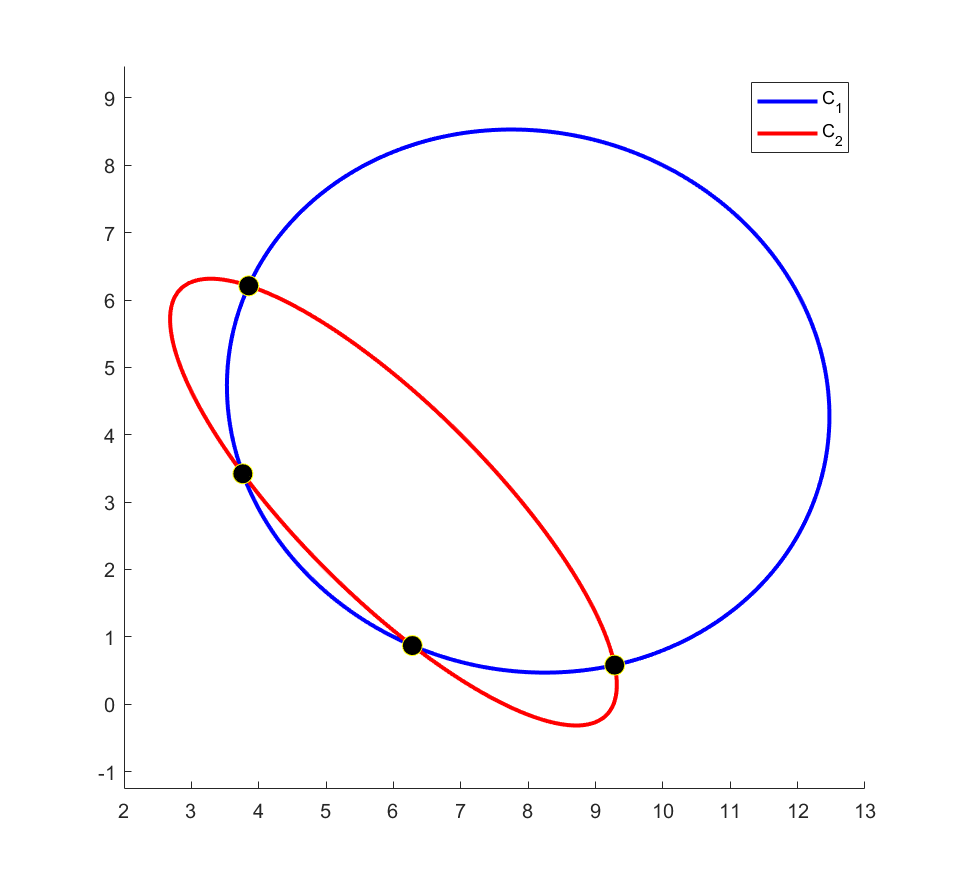}
    \caption{Two conics intersecting in four real points.}
    \label{fig:parabola-1}
\end{figure}

Consider the geometry of Fig.~\ref{fig:parabola-1}, where the two conics are associated with the matrices 
\[
    \bC_1 \propto \begin{bmatrix}
       65&4&-538\\
           4&80&-392\\
        -538&-392&4772
    \end{bmatrix}\qquad 
\bC_2 \propto \begin{bmatrix}
    11&9&-93\\
     9&11&-87\\
   -93&-87&779
\end{bmatrix}
\]

\subsubsection*{Intersection with the canonical parameterization}
From the two conic matrices given above, we can select three arbitrary points using Eq.~\eqref{eq:points-on-conic}
\[
\bp_1 \propto \begin{bmatrix}
 0 \\
          4.9 +     5.9699i\\
            1 
 \end{bmatrix} \qquad \bp_2 \propto \begin{bmatrix}
   0 \\
          4.9 -     5.9699i\\
            1\\
 \end{bmatrix} \qquad \bp_3 \propto \begin{bmatrix}
      8.2769 +     2.2154i\\
           0\\
            1 
 \end{bmatrix}
\]
From \(\bp_1\) and \(\bp_2\) we can obtain $\bp_0$ from Eq.~\eqref{eq:p0}
\[
\bp_0 \propto \begin{bmatrix}
    5.5\\
        4.625\\
            1
\end{bmatrix}
\]
We can now solve for \(\lambda_{0,1,2}\) taking the inverse of the matrix on the left-hand side of Eq.~\eqref{eq:lambdas}
\[
\begin{bmatrix}
    \lambda_0\\
    \lambda_1\\
    \lambda_2
\end{bmatrix} = \begin{bmatrix}
    \bp_0 & \bp_1 & \bp_2
\end{bmatrix}^{-1}\bp_3 =\begin{bmatrix}  1.5049 +     0.4028i\\
     -0.24317 +    0.17433i\\
     -0.26172 -    0.57713i
\end{bmatrix}
\]
so that, with Eq.~\eqref{eq:H}, we obtain
\[
    \bH \propto \begin{bmatrix}
       8.2769 +2.2154i&0 &0 \\
       6.9601 +1.8629i&-2.2323 -    0.59749i&-4.7279 -1.2655i\\
       1.5049 +0.4028i&-0.24317 +0.17433i &    -0.26172 -0.57713i
    \end{bmatrix}
\]
Using this transformation, the conic matrices become (Eq.~\eqref{eq:c1p-parabola} and Eq.~\eqref{eq:c2Generic})
\[
\bC_1' \propto \begin{bmatrix}
2  & 0 & 0\\
0 & 0 &-1\\
0 &  -1  &0
\end{bmatrix}\qquad 
\bC_2'\propto \begin{bmatrix}
   -0.4800 - 0.2768i&0.1712 - 0.4623i&-0.6658 + 0.8042i\\
   0.1712 - 0.4623i&-0.3930 + 0.0573i&0.9569 + 0.5518i\\
  -0.6658 + 0.8042i&0.9569 + 0.5518i&-0.6607 - 1.6544i\\
\end{bmatrix}
\]
With Eq.~\eqref{eq:quartic}, we can now determine the quartic that we need to solve
\[
(-0.3930 + 0.0573i) x^4+(0.3423 - 0.9245i)x^3+(1.4338 + 0.8268i)x^2+(-1.3315 + 1.6084i)x-0.6607 - 1.6544i =0
\]
whose solutions are the four numbers
\[
x_{a}' = 1.4431 + 0.1883i\quad x_b' = -1.4305 - 0.2675i \quad  x_c' =1.3020 - 0.6502i  \quad x_d'=-0.1259 - 1.4499i
\]
The \(y\)-components can be calculated using Eq.~\eqref{eq:parabola}
\[
y_a'= 2.0471 + 0.5433i\quad y_b'=1.9749 + 0.7652i\quad y_c' = 1.2725 - 1.6931i\quad y_d' = -2.0863 + 0.3651i
\]
Finally, with Eq.~\eqref{eq:quartic-s} 
\[
\bs_a\propto \begin{bmatrix}
9.2839\\
    0.5803\\
    1
   \end{bmatrix}\quad \bs_b \propto \begin{bmatrix}
        3.8515\\
    6.2106\\
    1
\end{bmatrix}\quad \bs_c \propto \begin{bmatrix}
   6.2804\\
    0.8705\\
    1
\end{bmatrix}\quad \bs_d \propto \begin{bmatrix}
 3.7641\\
    3.4209\\
    1
\end{bmatrix}
\]
where the first two components are the coordinates of the intersection points in \(\mathbb{R}^2\).

\subsubsection*{Intersection with the self-polar triangle}
Using the same conics reported above, we may calculate
\[
\bC_2^*\bC_1\propto\begin{bmatrix}-4.0625&-0.2500&12.4946  \\
    0.3614&-3.3370&1.5598\\
   -0.3696&-0.3478&1
\end{bmatrix}
\]
whose eigenvectors are obtained solving Eq.~\eqref{eq:self-eig}
\[
\bp_0 \propto \begin{bmatrix}
     -0.9263\\
   -0.2811\\
   -0.2509
\end{bmatrix}\qquad \bp_1\propto \begin{bmatrix}
    0.8987\\
    0.4075\\
    0.1624
\end{bmatrix}\qquad \bp_2 \propto \begin{bmatrix}
  -0.7222\\
    0.6917\\
   -0.0056
\end{bmatrix}
\]
Using
\[
\bp_3\propto \begin{bmatrix}
    1\\
    1\\
    1
\end{bmatrix}
\]
we solve the system of Eq.~\eqref{eq:lambdas}
\[
\begin{bmatrix}
    \lambda_0\\
    \lambda_1\\
    \lambda_2
\end{bmatrix} = \begin{bmatrix}
    -7.4823\\
   -5.3493\\
    1.5559
\end{bmatrix}
\]
so that, using Eq.~\eqref{eq:H}
\[
\bH \propto \begin{bmatrix}
    6.9308&-4.8072&-1.1236\\
    2.1036&-2.1798&1.0762\\
    1.8772&-0.8685&-0.0087
\end{bmatrix}
\]
As expected, the conic matrices are transformed to the diagonal matrices (see Eq.~\eqref{eq:c1p} and Eq.~\eqref{eq:c2p})
\[
\bC_1'\propto \begin{bmatrix}
     3.3138  & 0 &   0\\
   0 &  -0.4110  &  0\\
    0   & 0  &  0.1622\\
\end{bmatrix}\qquad \bC_2\propto \begin{bmatrix}
     0.4776 &  0  &0\\
         0  & -0.0233    &0\\
   0  &  0   & 0.0047
\end{bmatrix}
\]
and we can easily find the solution to the system of equations composed of Eq.~\eqref{eq:self-polar1} and Eq.~\eqref{eq:self-polar2}
\[
3.3138 x'^2 -0.4110 y'^2+0.1622=0
\]
\[
0.4776x'^2-0.0233y'^2+0..0047=0
\]
that are given by
\[
x' = \pm  0.1241\qquad y' = \pm 0.7203
\]
From these values we can write down the coordinates of the points of intersection in the new frame (Eq.~\eqref{eq:self-s_p}:
\[
\bs_a' \propto \begin{bmatrix}
    0.1241\\
    0.7203\\
    1
\end{bmatrix}\quad\bs_b' \propto \begin{bmatrix}
    -0.1241\\
    -0.7203\\
    1
\end{bmatrix}\quad\bs_c' \propto \begin{bmatrix}
    0.1241\\
    -0.7203\\
    1
\end{bmatrix}\quad\bs_d' \propto \begin{bmatrix}
    -0.1241\\
    0.7203\\
    1
\end{bmatrix}
\]
that give, in the original coordinates (Eq.~\eqref{eq:self-s})
\[
    \bs_a \propto \begin{bmatrix}
        9.2839\\
    0.5803\\
    1
    \end{bmatrix}\quad \bs_b\propto \begin{bmatrix}
        3.8515\\
    6.2106\\
    1
    \end{bmatrix}\quad \bs_c\propto \begin{bmatrix}
        3.7641\\
    3.4209\\
    1
    \end{bmatrix}\quad \bs_d\propto \begin{bmatrix}
         6.2804\\
    0.8705\\
    1
    \end{bmatrix}
\]

We can see that the points obtained with this approach coincide with the points obtained solving the quartic equation, and direct substitution into the original equations shows that they all belong to the two conics. 

\subsection{Two touching conics}
\begin{figure}[ht!]
    \centering
\includegraphics[width=0.5\textwidth]{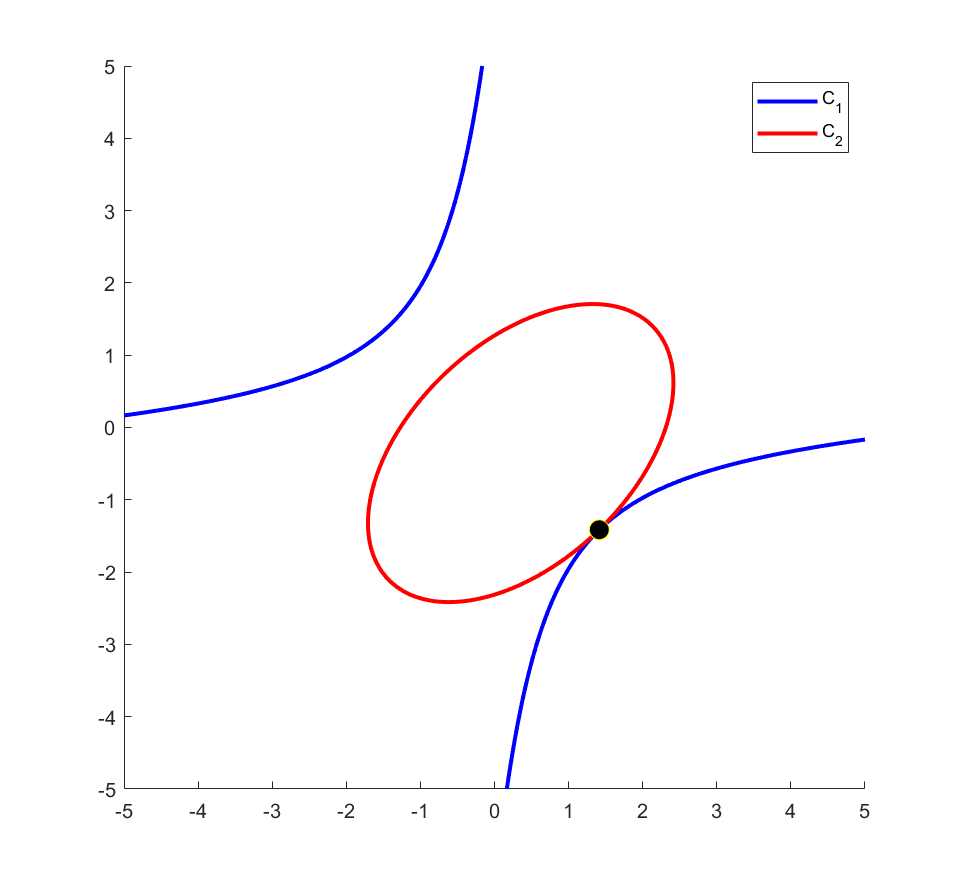}
    \caption{Two conics tangent at one real point. }
    \label{fig:tangent}
\end{figure}

Consider the two conics of Fig.~\ref{fig:tangent}. Their matrices are given by
\[
\bC_1 \propto \begin{bmatrix}
    -8 & 72 & 0\\
    672 & -8 & 0\\
    0 & 0 & 320
\end{bmatrix}\qquad \bC_2 \propto \begin{bmatrix}
    68 & -32 & -25\sqrt{2}\\
    -32 & 68 & 25\sqrt{2}\\
    -25 \sqrt{2} & 25\sqrt{2} & -200
\end{bmatrix}
\]

In this case, the intersection performed with the quartic equation does not provide any special insights, as it applies as in the previous example. Instead, we will analyze how to deal with this configuration using the approach based on the self-polar triangle.

In this case, two of the vertices of the self-polar triangle coincide. We let this point be
\[
\bp_1 \propto \begin{bmatrix}
    1.4142\\
    -1.4142\\
    1
\end{bmatrix}
\]
and select \(\bp_2\) and \(\bp_3\) on the conic using any two points from Eq.~\eqref{eq:points-on-conic}:
\[
\bp_2 \propto \begin{bmatrix}
    0\\
    6.3246\\
    1.0000
    \end{bmatrix}\qquad \bp_3 \propto
    \begin{bmatrix}-6.3246\\
         0\\
    1.0000
    \end{bmatrix}
\]
Now we evaluate \(\bp_0\) using Eq.~\eqref{eq:p0}:
\[
\bp_0\propto \begin{bmatrix}
    -1.1441\\
   -3.9725\\
    1
\end{bmatrix}
\]
and proceed as usual to determine the matrix \(\bH\) with Eq.~\eqref{eq:H}:
\[
    \bH \propto \begin{bmatrix}
        -2.8284& -3.4961 &       0\\
   -9.8207& 3.4961& 6.3246\\
    2.4721&-2.4721&1
    \end{bmatrix}
\]
from which
\[
\bC_1 \propto \begin{bmatrix}5120   &0 &0\\
0   &0 &    -2560\\
    0&-2560&0\\ \end{bmatrix}\qquad \bC_2 \propto \begin{bmatrix}
        2880  &0&-3840\\
  0&0&      2400\\
        -3840&2400&2967.2
    \end{bmatrix}
\]
The quadratic equation that we need to solve is given by Eq.~\eqref{eq:quadratic}
\[
7680 x'^2  -7680x'+2967.2=0
\]
from which, using Eq.~\eqref{eq:quadratic-sol} and Eq.~\eqref{eq:parabola}
\[
x_{a,b}'=0.5\pm0.3693i,\quad y_{a,b}'= 0.11364\pm  0.3693i
\]
The four points of intersection are finally given by
\[
\bs_{s,b}\propto  \begin{bmatrix} 
-0.92655 \mp     1.1945i \\
0.92655 \mp    1.1945i\\
1
\end{bmatrix}\qquad \bs_{c,d} \propto \bp_1
\]

\section{Conclusions}
In this work, we have shown that the intersection between two conics can be easily accomplished through a convenient change of coordinates. The two solutions presented allow the analyst to solve the intersection problem by transformation into canonical coordinates (requiring the solution to a quartic equation) or a set of coordinates defined by the conics' common self-polar triangle (requiring the solution to a quadratic equation). Both methods are shown to be effective in two numerical examples, though the self-polar triangle approach provides a simpler algorithm.

\newpage

\begin{appendices}

\section{Pseudo-codes for conic intersection}\label{sec:pseudocodes}

\begin{algorithm}[ht!]
\caption{Intersection with canonical representation}\label{alg:method2}
\begin{algorithmic}[1]
\State \textbf{given} conics \(\bC_1\) and \(\bC_2\)
\State compute the coordinates of three points lying on it \(\bC_1\): \(\bp_1\), \(\bp_2\) and \(\bp_3\) \Comment{Eq.~\eqref{eq:points-on-conic}}
\State compute \(\bp_0 \propto \left(\bC_1\bp_1\right)\times \left(\bC_1\bp_2\right)\)  \Comment{Eq.~\eqref{eq:p0}}
\State set \(\lambda_3 = 1\) and compute \(\lambda_{1,2,3}\) \Comment{Eq.~\eqref{eq:lambdas}}
\State compute \(\bH\) \Comment{Eq.~\eqref{eq:H}}
\State compute \(\bC_2'\) \Comment{Eq.~\eqref{eq:c2Generic}}
\State compute \(x'_{a,b,c,d}\) solving the quartic of Eq.~\eqref{eq:quartic} \Comment{see Appendix~\ref{app:quarticsolution}}
\State compute \(y'_{a,b,c,d}\) \Comment{Eq.~\eqref{eq:parabola}}
\State form \(\bs_j' \propto [x_j\,y_j\,1]^T\), \(j=a,b,c,d\) \Comment{Eq.~\eqref{eq:quartic-sp}}
\State compute \(\bs_j\propto\bH\bs_j'\), \(j=a,b,c,d\) \Comment{Eq.~\eqref{eq:quartic-s}}
\State \textbf{output} intersection points \(\bs_j\), \(j=a,b,c,d\)
\end{algorithmic}
\end{algorithm}

\begin{algorithm}[ht!]
\caption{Intersection with self-polar triangle}\label{alg:cap}
\begin{algorithmic}[1]
\State \textbf{given} conics \(\bC_1\) and \(\bC_2\)
\State compute the eigenvectors \(\bp_{0,1,2}\) of \(\bC_2^{*}\bC_1\)  \Comment{Eq.~\eqref{eq:self-eig}}
\If{\(\bp_l^T\bC_1\bp_l=\bp_l^T\bC_2\bp_l=\bp_m^T\bC_1\bp_m=\bp_m^T\bC_2\bp_m=0\) for some \(l\neq m\)}
\State set \(\bp_1=\bp_l\), select \(\bp_2\), \(\bp_3\) on the conic \Comment{Eq.~\eqref{eq:points-on-conic}}
\State compute \(\bp_0\propto (\bC_1\bp_1)\times (\bC_1\bp_2)\) \Comment{Eq.~\eqref{eq:p0}}
\State set \(\lambda_3=1\), compute \(\lambda_{0,1,2}\)\Comment{Eq.~\eqref{eq:lambdas}}
\State compute \(\bH\) \Comment{Eq.~\eqref{eq:H}}
\State compute \(\bC_2'\propto \bH^T\bC_2\bH\)
\State compute \(x_{a,b}'\) \Comment{Eq.~\eqref{eq:quadratic-sol}}
\State compute \(y_{a,b}'=x_{a,b}'^2\) 
\State set \(\bs_a'^T\propto [x_a\,y_a\,1]\), \(\bs_b'^T\propto [x_b\,y_b\,1]\)
\State evaluate \(\bs_{a,b}\propto \bH\bs_{a,b}'\)
\State \textbf{output} intersection points \(\bp_1, \bs_{a,b}\)
\EndIf 
\State set \(\bp_{3}\propto[1\,1\,1]^T,\,\lambda_3=1\)
\State compute \(\lambda_{0,1,2}\) \Comment{Eq.~\eqref{eq:lambdas}}
\State compute \(\bH\) \Comment{Eq.~\eqref{eq:H}}
\State compute \(\bC_1'\) and \(\bC_2'\) \Comment{Eq.~\eqref{eq:c1p},\ref{eq:c2p}}
\State compute \(x'^2\) and \(y'^2\) \Comment{Eq.~\eqref{eq:self-polarX_sq}, Eq.~\eqref{eq:self-polarY_sq_v2}}
\State compute \(\bs_j'\), \(j=a,b,c,d\) \Comment{Eq.~\eqref{eq:self-s_p}}
\State compute \(\bs_j\propto\bH\bs_j'\), \(j=a,b,c,d\) \Comment{Eq.~\eqref{eq:self-s}}
\State \textbf{output} intersection points \(\bs_j\), \(j=a,b,c,d\)
\end{algorithmic}
\end{algorithm}

\newpage
\section{Solution to a quartic equation}\label{app:quarticsolution}
Consider the quartic equation
\begin{equation}
    ax^4+bx^3+cx^2+dx+e=0
\end{equation}
where we assume \(a\neq 0\). Introduce the parameters
\begin{equation}
    p = \frac{8ac-3b^2}{8a^2}, \qquad S = \frac{8a^2d-4abc+b^3}{8a^3}
\end{equation}
\begin{equation}
    q = 12ae-3bd+c^2, \qquad 
    s = 27ad^2-72ace+27b^2e-9bcd+2c^3
\end{equation}
\begin{equation}
    \Delta_0 = \sqrt[3]{\frac{s+\sqrt{s^2-4q^3}}{2}},\qquad
    Q = \frac{1}{2}\sqrt{-\frac{2}{3}p+\frac{1}{3a}\left(\Delta_0+\frac{q}{\Delta_0}\right)}
\end{equation}
Note that there are three possible values for \(\Delta_0\). If \(Q\simeq 0\) for the chosen \(\Delta_0\), then it should be changed to 
\begin{equation}
    \Delta_0 = \Delta_0 e^{2\pi n i / 3}\qquad n\in \mathbb{Z}
\end{equation}
and the corresponding \(Q\) updated with this new value. 
The solutions to the equation are finally given by
\begin{equation}
    x_{1,2}=-\frac{b}{4a}-Q\pm \frac{1}{2}\sqrt{-4Q^2-2p+\frac{S}{Q}}
\end{equation}
\begin{equation}
    x_{3,4}=-\frac{b}{4a}+Q\pm \frac{1}{2}\sqrt{-4Q^2-2p-\frac{S}{Q}}
\end{equation}

\end{appendices}

\printbibliography

\end{document}